\documentclass[12pt,a4paper,twoside]{article}

\newcommand{\pageformat}[6]{\setlength{\hoffset}{-1in}
                  \setlength{\voffset}{-1in}
                  \addtolength{\hoffset}{#5}
                            \addtolength{\voffset}{#6}
                            \setlength{\oddsidemargin}{#1}
                            \setlength{\evensidemargin}{#2}
                            \setlength{\textwidth}{\paperwidth}
                  \addtolength{\textwidth}{-\oddsidemargin}
                  \addtolength{\textwidth}{-\evensidemargin}
                  \addtolength{\textwidth}{-\marginparsep}
                  \addtolength{\textwidth}{-\marginparwidth}
                            \setlength{\topmargin}{#3}
                            \setlength{\textheight}{\paperheight}
                  \addtolength{\textheight}{-\topmargin}
                  \addtolength{\textheight}{-\headheight}
                  \addtolength{\textheight}{-\headsep}
                  \addtolength{\textheight}{-\footskip}
                  \addtolength{\textheight}{-#4}}
\pageformat{2cm}{3cm}{24mm}{24mm}{1pt}{0pt}

\usepackage{ifthen}
\newboolean{article}
    \setboolean{article}{true}
\newboolean{report}
\newboolean{book}
\newboolean{letter}
\newboolean{german}
\newboolean{italian}
\newboolean{nobaselinestretch}
\newboolean{nosectionappendix}
\newboolean{oldtoc}
\newboolean{nosectionequn}
\newboolean{notheorem}

\ifthenelse{\boolean{german}}{
    \usepackage{german}}{}

\usepackage[latin1]{inputenc}

\usepackage{amsmath}
\usepackage{amssymb}
\usepackage[mathscr]{eucal}

\ifthenelse{\boolean{notheorem}}{}{
    \usepackage{theorem}}

\ifthenelse{\boolean{nobaselinestretch}}{}{
    \renewcommand{\baselinestretch}{1.25}}

\newenvironment{env}[2]{\begin{#1}#2\end{#1}}{}
    \newcommand{\beq}[1]{\begin{env}{equation}{#1}}
    \newcommand{\beqn}[1]{\begin{env}{equation*}{#1}}
    \newcommand{\bal}[1]{\begin{env}{align}{#1}}
    \newcommand{\baln}[1]{\begin{env}{align*}{#1}}
    \newcommand{\bga}[1]{\begin{env}{gather}{#1}}
    \newcommand{\bgan}[1]{\begin{env}{gather*}{#1}}
    \newcommand{\bflal}[1]{\begin{env}{flalign}{#1}}
    \newcommand{\bflaln}[1]{\begin{env}{flalign*}{#1}}
    \newcommand{\bmu}[1]{\begin{env}{multline}{#1}}
    \newcommand{\bmun}[1]{\begin{env}{multline*}{#1}}
    \newcommand{\bsp}[1]{\begin{env}{split}{#1}}

    \newcommand{\eeq}{\end{env}}
    \newcommand{\eeqn}{\end{env}}
    \newcommand{\eal}{\end{env}}
    \newcommand{\ealn}{\end{env}}
    \newcommand{\ega}{\end{env}}
    \newcommand{\egan}{\end{env}}
    \newcommand{\eflal}{\end{env}}
    \newcommand{\eflaln}{\end{env}}
    \newcommand{\emu}{\end{env}}
    \newcommand{\emun}{\end{env}}
    \newcommand{\esp}{\end{env}}

\newcommand{\lf}{\vspace{2ex}}

\renewcommand{\bf}[1]{\textbf{#1}}
\renewcommand{\it}[1]{\textit{#1}}

\renewcommand{\sf}[1]{\textsf{#1}}

\renewcommand{\tt}[1]{\texttt{#1}}
\newcommand{\hl}[1]{\bf{\it{#1}}}

\newcommand{\mbf}[1]{\mathbf{#1}}
\newcommand{\msf}[1]{\text{\small$\sf{#1}$}}

\newcommand{\cmc}[1]{\mathcal{#1}}
\newcommand{\eus}[1]{\mathscr{#1}}
\newcommand{\euf}[1]{\mathfrak{#1}}
\newcommand{\bb}[1]{\mathbb{#1}}

\newcommand{\mfootnotesize}[1]{{\setlength{\arraycolsep}{.5ex}\text{\footnotesize$#1$}}}
\newcommand{\mscriptsize}[1]{{\setlength{\arraycolsep}{.3ex}\text{\scriptsize$#1$}}}

\newcommand{\nbd}[1]{$#1$\nobreakdash--}
\newcommand{\ol}[1]{\overline{#1}}

\newcommand{\wt}[1]{\widetilde{#1}}

\newcommand{\vt}{\vartheta}

\newcommand{\AB}[1]{\langle#1\rangle}

\newcommand{\CB}[1]{\{#1\}}
\newcommand{\bCB}[1]{\bigl\{#1\bigr\}}
\newcommand{\BCB}[1]{\Bigl\{#1\Bigr\}}
\newcommand{\SB}[1]{[#1]}

\newcommand{\Matrix}[1]{\begin{pmatrix}#1\end{pmatrix}}

\newcommand{\fMatrix}[1]{\mfootnotesize{\Matrix{#1}}}
\newcommand{\sMatrix}[1]{\mscriptsize{\Matrix{#1}}}

\newcommand{\sbars}[1]{\:\bar{#1}^s\:}
\newcommand{\sodots}{\sbars{\odot}}
\newcommand{\set}[2][]{
    \ifthenelse{\equal{#1}{}}{
        \CB{#2}}{
        \CB{#1~|~#2}}}
\newcommand{\bset}[2][]{
    \ifthenelse{\equal{#1}{}}{
        \bCB{#2}}{
        \bCB{#1~|~#2}}}
\newcommand{\Bset}[2][]{
    \ifthenelse{\equal{#1}{}}{
        \BCB{#2}}{
        \BCB{#1~\big|~#2}}}
\newcommand{\zero}{\CB{0}}

\DeclareMathOperator{\ls}{\normalfont\msf{span}}
\DeclareMathOperator{\cls}{\ol{\ls}}

\DeclareMathOperator{\id}{\normalfont\msf{id}}

\newcommand{\C}{\bb{C}}

\newcommand{\E}{\bb{E}}

\newcommand{\cA}{\cmc{A}}
\newcommand{\cB}{\cmc{B}}
\newcommand{\cC}{\cmc{C}}
\newcommand{\cD}{\cmc{D}}

\newcommand{\cM}{\cmc{M}}

\newcommand{\sB}{\eus{B}}

\newcommand{\sF}{\eus{F}}

\newcommand{\sK}{\eus{K}}

\newcommand{\el}{\euf{l}}

\newcommand{\er}{\euf{r}}

\newcommand{\eC}{\euf{C}}

\newcommand{\U}{\mbf{1}}

\newcommand{\f}{\text{\scriptsize$\sF$}}

\ifthenelse{\boolean{nosectionequn}}{}{
    \numberwithin{equation}{section}
    }

\ifthenelse{\boolean{article}\or\boolean{letter}\or\boolean{nosectionequn}}{
    \setboolean{nosectionappendix}{true}}{}
\ifthenelse{\boolean{nosectionappendix}}{}{
    \renewcommand{\appendix}{
        \chapter*{\appendixname}
        \addcontentsline{toc}{chapter}{\appendixname}
        \renewcommand{\thesection}{\Alph{section}}
        \setcounter{section}{0}}}
   
\ifthenelse{\boolean{report}\or\boolean{book}}{
    }{}

\ifthenelse{\boolean{notheorem}}{}{
        \newcommand{\mnname}{Mathematical note.}
        \newcommand{\enname}{End of the note.}
        \newcommand{\definame}{Definition.}
        \newcommand{\propname}{Proposition.}
        \newcommand{\lemname}{Lemma.}
        \newcommand{\exname}{Example.}
        \newcommand{\exername}{Exercise.}
        \newcommand{\remname}{Remark.}
        \newcommand{\obname}{Observation.}
        \newcommand{\thmname}{Theorem.}
        \newcommand{\corname}{Corollary.}
        \newcommand{\proofname}{Proof.}
    \ifthenelse{\boolean{german}}{
        \renewcommand{\mnname}{Mathematische Notiz.}
        \renewcommand{\enname}{Ende der Notiz.}
        \renewcommand{\exname}{Beispiel.}
        \renewcommand{\exername}{Übung.}
        \renewcommand{\remname}{Bemerkung.}
        \renewcommand{\obname}{Beobachtung.}
        \renewcommand{\thmname}{Satz.}
        \renewcommand{\corname}{Korollar.}
        \renewcommand{\proofname}{Beweis.}}{}
    \ifthenelse{\boolean{italian}}{
        \renewcommand{\mnname}{Nota matematica.}
        \renewcommand{\enname}{Fina della nota.}
        \renewcommand{\definame}{Definizione.}
        \renewcommand{\propname}{Proposizione.}
        \renewcommand{\exname}{Esempio.}
        \renewcommand{\exername}{Esercizio.}
        \renewcommand{\remname}{Nota.}
        \renewcommand{\obname}{Osservazione.}
        \renewcommand{\thmname}{Teorema.}
        \renewcommand{\corname}{Corollario.}
        \renewcommand{\proofname}{Dimostrazione.}

       \renewcommand{\appendixname}{Appendice}

       }{}
    \theoremheaderfont{\normalfont\bfseries}
    \theoremstyle{change}
        \theorembodyfont{\rmfamily}
            \newtheorem{emp}{}[section]
                \newcommand{\bemp}[1][]{
                    \begin{emp}\hskip-\labelsep\bf{#1}\hskip\labelsep}
                \newcommand{\eemp}{\end{emp}}
\newtheorem{itemp}[emp]{}
                \newcommand{\bitemp}[1][]{
                    \begin{itemp}\hskip-\labelsep\bf{#1}\hskip\labelsep\normalfont\itshape}
                \newcommand{\eitemp}{\end{itemp}}
            \newtheorem{mn}[emp]{\mnname}
                \newcommand{\bnm}{\begin{mn}~\begin{quotation}\renewcommand{\baselinestretch}{1}\small\noindent\ignorespaces}
                \newcommand{\enm}{\end{quotation}\hfill\bf{\enname}\end{mn}}
            \newtheorem{ex}[emp]{\exname}
                \newcommand{\bex}{\begin{ex}}
                \newcommand{\eex}{\end{ex}}
            \newtheorem{exer}[emp]{\exername}
                \newcommand{\bexer}{\begin{exer}}
                \newcommand{\eexer}{\end{exer}}
            \newtheorem{defi}[emp]{\definame}
                \newcommand{\bdefi}{\begin{defi}}
                \newcommand{\edefi}{\end{defi}}
            \newtheorem{rem}[emp]{\remname}
                \newcommand{\brem}{\begin{rem}}
                \newcommand{\erem}{\end{rem}}
            \newtheorem{ob}[emp]{\obname}
                \newcommand{\bob}{\begin{ob}}
                \newcommand{\eob}{\end{ob}}
        \theorembodyfont{\normalfont\itshape}
            \newtheorem{thm}[emp]{\thmname}
                \newcommand{\bthm}{\begin{thm}}
                \newcommand{\ethm}{\end{thm}}
            \newtheorem{prop}[emp]{\propname}
                \newcommand{\bprop}{\begin{prop}}
                \newcommand{\eprop}{\end{prop}}
            \newtheorem{cor}[emp]{\corname}
                \newcommand{\bcor}{\begin{cor}}
                \newcommand{\ecor}{\end{cor}}
            \newtheorem{lem}[emp]{\lemname}
                \newcommand{\blem}{\begin{lem}}
                \newcommand{\elem}{\end{lem}}
\newenvironment{empn}[1]{\lf\noindent\bf{#1}\ignorespaces\hskip\labelsep}{\lf}
		\newcommand{\bempn}[1]{\begin{empn}{#1}}
		\newcommand{\eempn}{\end{empn}}
		\newcommand{\bitempn}[1]{\begin{empn}{#1}\normalfont\itshape}
		\newcommand{\eitempn}{\end{empn}}
                \newcommand{\bnmn}{\begin{empn}{\mnname}~\begin{quotation}\renewcommand{\baselinestretch}{1}\small\noindent\ignorespaces}
                \newcommand{\enmn}{\end{quotation}\hfill\bf{\enname}\end{empn}}
		\newcommand{\bexn}{\begin{empn}{\exname}}
		\newcommand{\eexn}{\end{empn}}
		\newcommand{\bexern}{\begin{empn}{\exername}}
		\newcommand{\eexern}{\end{empn}}
		\newcommand{\bdefin}{\begin{empn}{\definame}}
		\newcommand{\edefin}{\end{empn}}
		\newcommand{\bremn}{\begin{empn}{\remname}}
		\newcommand{\eremn}{\end{empn}}
		\newcommand{\bobn}{\begin{empn}{\obname}}
		\newcommand{\eobn}{\end{empn}}

\newcommand{\qedsymbol}{~\rule[-0.35mm]{2mm}{2mm}}
    \newcounter{proof}[emp]
    \newenvironment{Proof}[1]{
        \vspace{1ex}
        \renewcommand{\item}[1][\stepcounter{proof}(\roman{proof})]%
            {##1\hskip\labelsep}
        \noindent\textsc{#1\hskip\labelsep}}{
        \nolinebreak\qedsymbol}
    \newcommand{\proof}[1][\proofname]{
        \begin{Proof}{#1}\ignorespaces}
    \newcommand{\qed}{\end{Proof}}
    \newcommand{\noqed}{
        \renewcommand{\qedsymbol}{}
        \end{Proof}}}
    \ifthenelse{\boolean{italian}}{
        \renewcommand{\proofname}{Dimostrazione.}}{}

\usepackage{txfonts}

\usepackage[matrix,arrow,curve]{xy}

\newcommand{\cvN}{\euf{cvN}}

\begin{document}

\title{Commutants of von Neumann Correspondences\\and Duality of\\Eilenberg-Watts Theorems by Rieffel and by Blecher}
\author{}
\author{
~\\
Michael Skeide\thanks{This work is supported by research fonds of the Department S.E.G.e S.\ of University of Molise.}\\\\
{\small\itshape Dipartimento S.E.G.e S.}\\
{\small\itshape Università degli Studi del Molise}\\
{\small\itshape Via de Sanctis}\\
{\small\itshape 86100 Campobasso, Italy}\\
{\small{\itshape E-mail: \tt{skeide@math.tu-cottbus.de}}}\\
{\small{\itshape Homepage: \tt{http://www.math.tu-cottbus.de/INSTITUT/lswas/\_skeide.html}}}\\
\\
}
\date{January 2005}

{
\renewcommand{\baselinestretch}{1}
\maketitle

\vfill

\begin{abstract}
\noindent
The category of von Neumann correspondences from $\cB$ to $\cC$ (or von Neumann \nbd{\cB}\nbd{\cC}mod\-ules) is \it{dual} to the category of von Neumann correspondences from $\cC'$ to $\cB'$ via a functor that generalizes naturally the functor that sends a von Neumann algebra to its commutant and back. We show that under this duality, called \it{commutant}, Rieffel's \it{Eilenberg-Watts theorem} (on functors between the categories of representations of two von Neumann algebras) switches into Blecher's \it{Eilenberg-Watts theorem} (on functors between the categories of von Neumann modules over two von Neumann algebras) and back.
\end{abstract}

\vfill

}




\section{Introduction}

In algebra the \it{Eilenberg-Watts theorem} (Eilenberg \cite{Eil60} and Watts \cite{Wat60}) states that every functor (fulfilling certain assumptions) between categories of (right or left) modules over algebras $\cB$ and $\cC$ is given (up to natural equivalence) by tensoring with a (unique up to isomorphism) bimodule (from the right or from the left). As a consequence, such a functor is an equivalence, if and only if it is implemented by a bimodule that has an \it{inverse} under taking tensor products. This is the famous \it{Morita theorem} (Morita \cite{Mor58}) that answers the question under which conditions two algebras have equivalent representation theories.

In Hilbert module theory the symmetry between left and right actions is broken by existence of an inner product with values in that algebra which acts on the right.%
\footnote{
We are not speaking about the possibility to model the complete theory in terms of left Hilbert modules. A left Hilbert module would be a left module over a \nbd{C^*}algebra with an inner product assuming values in that \nbd{C^*}algebra and that inner product should be left linear with respect to its left argument. Our Hilbert modules will always be \hl{right} Hilbert modules with an inner product that is \hl{right} linear in its \hl{right} variable. The whole theory based on left modules will be symmetric and will have formulae equally elegant as those in the theory for right modules, only if we would decide (as, in fact, some algebraists do) to write all functions on the right of their argument. (An instance of this fact can be seen already in the notations for rank-one operators in Hilbert space theory, when people insist in having inner products linear in the first and not in the second variable.)

Actually the symmetry between left and right modules is already broken in the very moment when we choose to write functions on the left of their arguments. In fact, the condition to be a module map reads for right modules $T(xb)=(Tx)b$, a simple associativity where all parts of the formula stay in order. For left modules the condition reads $T(bx)=b(Tx)$. Here something has to be commuted. Also, in both parts the simple change of brackets would lead to expressions $Tb$ and $bT$ that, \it{a priori}, are not even defined.
}
Consequently, in Hilbert module theory there are two \it{Eilenberg-Watts theorems}.

The first one, Rieffel's \it{Eilenberg-Watts theorem} \cite{Rie74a}, actually concerns representations of von Neumann algebras on Hilbert spaces.%
\footnote{
Some authors, in fact, understand by a Hilbert module simply a representation of a \nbd{C^*}algebra on a Hilbert space, thus, turning the Hilbert space into a left module. This is definitely not the terminology we are using. For us a Hilbert module over a \nbd{\cB}algebra $\cB$ is a \hl{right} \nbd{\cB}module with a \nbd{\cB}valued inner product.
}
A Hilbert space is a von Neumann module over $\C$. True von Neumann modules appear only in  the statement of the theorem: Every \it{normal \nbd{*}functor} $\el$ from the category of representations%
\footnote{
Representations are assumed nondegenerate, unless stated otherwise explicitly.
}
of the von Neumann algebra $\cB$ to the category of representations of the von Neumann algebra $\cC$ is implemented (up to natural equivalence) as tensoring from the left with a von Neumann correspondence%
\footnote{
By a \hl{correpondence} from a \nbd{C^*}algebra $\cC$ to a \nbd{C^*}algebra $\cB$ we mean a Hilbert \nbd{\cB}module $E$ with a nondegenerate representation of $\cC$ by adjointable operators on $E$ that turns $E$ into a \nbd{\cC}\nbd{\cB}module. (The nondegeneracy condition for correspondences is essential and not shared by all authors!) In particular, every Hilbert \nbd{\cB}module is a correspondence from the \nbd{C^*}algebra $\sB^a(E)$ of all adjointable operators on $E$ to $\cB$. Von Neumann correspondences will be defined in detail later.
}
$E$ from $\cC$ to $\cB$ (uniquely determined by $\el$ up to isomorphism). If we think of a Hilbert space $H$ that carries a normal unital representation of $\cB$ as a left \nbd{\cB}module, then $\el$ carries the representation of $\cB$ on $H$ to the induced representation of $\cC$ on the left \nbd{\cC}module $\el(H)=E\odot H$ (acting as $c(x\odot h)=cx\odot h$), where $E\odot H$ is the tensor product (over $\cB$) of the \nbd{\cC}\nbd{\cB}correspondence $E$ and the \nbd{\cB}\nbd{\C}correspondence $H$.%
\footnote{
The tensor product of a correspondence $E$ from $\cA$ to $\cB$ and a correspondence $F$ from $\cB$ to $\cC$ is that unique \nbd{\cB}\nbd{\cC}correspondence $E\odot F$ that is generated as a Hilbert module by elementary tensors $x\odot y$ with inner products $\AB{x_1\odot y_1,x_2\odot y_2}=\AB{y_1,\AB{x_1,x_2}\odot y_2}$.
}
A morphism from $H_1$ to $H_2$, that is, a left \nbd{\cB}linear mapping $a$, is sent to the morphism $\el(a)=\id_E\odot a$.

The second one, Blecher's \it{Eilenberg-Watts theorems} \cite{Ble97}, regards Hilbert modules: Every \it{strict \nbd{*}functor} $\er$ from the category of Hilbert \nbd{\cB}modules to the category of Hilbert \nbd{\cC}modules is implemented (up to natural equivalence) as tensoring from the right with a correspondence $F$ from $\cB$ to $\cC$ (uniquely determined by $\er$ up to isomorphism). That is, $\er$ carries the Hilbert \nbd{\cB}module $E$ to $\er(E)=E\odot F$ and a morphism $a\in\sB^a(E_1,E_2)$ is sent to $\er(a)=a\odot\id_E$.

While Rieffel's theorem fails for \nbd{C^*}algebras, Blecher's theorem allows for an obvious modification to normal \nbd{*}functors between categories of von Neumann modules (or abstract \nbd{W^*}modules). The failure of Rieffel's theorem for \nbd{C^*}algebras can be seen, for instance, by Rieffel's argument based on \cite[Proposition 8.19]{Rie74a}: Commutative \nbd{C^*}algebras have equivalent categories of representations, if and only if their spectra are Borel isomorphic. However, such an equivalence can be implemented by a tensor product, if and only if these \nbd{C^*}algebras have homeomorphic spectra. That is, the corresponding \it{Morita theorem} (that would follow, if the \it{Eilenberg-Watts theorem} would hold) may fail.

In Section \ref{repsec} we report recent results from Muhly, Skeide and Solel \cite{MuSo05} about representations of $\sB^a(E)$ and how they can be used to furnish a proof of Blecher's \it{Eilenberg-Watts theorem} without using operator space technology. Then, we use a duality, the \it{commutant}, between von Neumann \nbd{\cB}\nbd{\cC}correpondences and von Neumann \nbd{\cC'}\nbd{\cB'}correspondences, to prove that the two \it{Eilenberg-Watts theorems} (in the von Neumann variant of Blecher's version) are duals of each other via the commutant. More precisely, the commutant $'$ takes von Neumann \nbd{\C}\nbd{\cB}correpondences (that is, von Neumann \nbd{\cB}modules) to von Neumann \nbd{\cB'}\nbd{\C}correpon\-dences (that is normal unital representations of $\cB'$), and back. Under this transition a functor $\el$ translates into a functor ~$\er=(')\circ\el\,\circ(')$, and back. Moreover, if $\er$ is implemented by the von Neumann \nbd{\cB}\nbd{\cC}correpondence $F$, then $\el$ is implemented by the von Neumann \nbd{\cC'}\nbd{\cB'}correspon\-dence $F'$, and conversly. Therefore, any proof of Blecher's \it{Eilenberg-Watts theorem} (von Neumann case), automatically, also proofs Rieffel's, and conversely.

The commutant was discussed in Skeide \cite{Ske03c} in the case $\cC=\cB$. It also appeared, independently, in Muhly and Solel \cite{MuSo04} and was generalized to different algebras in Muhly and Solel \cite{MuSo05}. The idea of \it{commutant} has already produced a bunch of new results \cite{Ske03p,GoSk03p,Ske05a,Ske04p} (in preparation \cite{Ske03p1}) and has put known results into a new perspective \cite{Dix54,Sti55,Arv69,AlHK78,Sau80,Sau83,Arv89,Arv89a,Bha96,Goh04,Hir04p,Ske00b,MuSo02}. Still, there is no end in sight. We would like to give an account of all those new perpectives. But this probably would be beyond the space available for this note. So we content ourselves to refer the reader to the discussions in \cite{GoSk03p,Ske05a,Ske04p} (most of \cite{Ske03p} we discuss in Section \ref{vNmsec}) and promis to do our best to finish soon \cite{Ske03p1}.

In Section \ref{vNmsec} we discuss von Neumann modules and how to think of them as representations of the commutant of the von Neumann algebra. In Section \ref{commsec} we discuss von Neumann correspondences and the construction of their commutants. Finally, in Section \ref{EWsec} we show that the commutant interchanges the two \it{Eilenberg-Watts theorems}. The categories of concrete von Neumann modules and concrete von Neumann correspondences, which we introduce exactly for that goal, appear here for the first time.

\lf\noindent
\bf{Convention.~}
A von Neumann algebra is always a concrete algebra of operators acting nondegenerately on a Hilbert space. If we intend a von Neumann algebra without a fixed identifying representation --- but we never ever do that --- then we should say \nbd{W^*}algebra. The same convention applies to von Neumann modules (as opposed to \nbd{W^*}modules) which we always consider as concrete modules of operators between two Hilbert spaces. Only in that way the functor \it{commutant} is a true bijection, that is, not only up to (canonical) isomorphism or (natural) equivalence.

\section{Representations of $\sB^a(E)$ and Blecher's Eilenberg-Watts theorem}\label{repsec}

Suppose $E$ is a Hilbert module over a \nbd{C^*}algebra $\cB$. Then $E^*$ may be viewed as a correspondence from $\cB$ to $\sB^a(E)$, the \hl{dual correspondence} of $E$, when equipped with the bimodule operations $bx^*a:=(a^*xb^*)^*$ and with inner product $\AB{x^*,y^*}=xy^*$. Here $x^*$ is interpreted as the operator $y\mapsto\AB{x,y}$ in $\sB^a(E,\cB)$ with adjoint $x\colon b\mapsto xb$. Consequently, $xy^*$ is the \hl{rank-one operator} $z\mapsto x\AB{y,z}$ in $\sB^a(E)$.

Actually, the inner product of $E^*$ takes values only in the pre-\nbd{C^*}algebra $\sF(E)$ of \hl{finite-rank operators} (an ideal in $\sB^a(E)$) that is spaned linearly by the rank-one operators. Therefore, $E^*$ may be viewed also as a correspondence from $\cB$ to the \nbd{C^*}algebra $\sK(E)$ of \hl{compact operators}, the closure of $\sF(E)$ in $\sB^a(E)$ and a closed ideal.

We observe that $E\odot E^*=\sK(E)$ via the canonical isomorphism $x\odot y^*\mapsto xy^*$, where $\sK(E)$ (like every \nbd{C^*}algebra) is viewed as the \hl{identity correspondence} from $\sK(E)$ to $\sK(E)$ with its natual bimodule structure and inner product $\AB{a_1,a_2}=a_1^*a_2$.

Suppose now that $\vt\colon\sB^a(E)\rightarrow\sB^a(F)$ is a strict unital homomorphism, where $F$ is a Hilbert module over a \nbd{C^*}algebra $\cC$. This means, in particular, that $F$ is a correspondence from $\sB^a(E)$ to $\cC$. Moreover the left action being \hl{strict}, means that already the action of $\sK(E)$ on $F$ is nondegenerate. (This is the only property we need. So we do not give a precise definition of a strict mapping.) In other words, by strictness $F$ is even a correspondence from $\sK(E)$ to $\cB$. Moreover, since $\vt(xy^*)z$ $(x,y\in E,z\in F)$ is a total subset of $F$, the action of $a\in\sB^a(E)$ on $F$ is already determined by its restriction to $\sK(E)$ as $\vt(a)\vt(xy^*)z=\vt((ax)y^*)z$ where $(ax)y^*\in\sK(E)$.

Nondegeneracy of the action of $\sK(E)$ on $F$ can be phrased as $\sK(E)\odot F=F$ via the canonical isomorphism $a\odot z\mapsto\vt(a)z$. Putting together this identification with $E\odot E^*=\sK(E)$, we obtain
\beqn{
F
~=~
\sK(E)\odot F
~=~
(E\odot E^*)\odot F
~=~
E\odot(E^*\odot F)
~=~
E\odot F_\vt,
}\eeqn
where we defined the correspondence $F_\vt$ from $\cB$ to $\cC$ as $F_\vt:=E^*\odot F$. The following theorem from Muhly, Skeide and Solel \cite{MSS03p} just fixes the isomorphism $F=E\odot F_\vt$ and identifies the action of $a\in\sB^a(E)$ on $E\odot F_\vt$ as the canonical one, that is, as amplification $\vt(a)=a\odot\id_{F_\vt}$.

\bitemp[Theorem \cite{MSS03p}.]\label{strirepthm}
Let $E$ be a Hilbert \nbd{\cB}module, let $F$ be a Hilbert \nbd{\cC}module and let $\vt\colon\sB^a(E)\rightarrow\sB^a(F)$ be a strict unital homomorphism. Then $F_\vt:=E^*\odot F$ is a correspondence from $\cB$ to $\cC$ and the formula
\beqn{
u(x_1\odot(x_2^*\odot y))
~:=~
\vt(x_1x_2^*)y
}\eeqn
defines a unitary
\beqn{
u
\colon
E\odot F_\vt
~\longrightarrow~
F
}\eeqn
such that
\beqn{
\vt(a)
~=~
u(a\odot\id_{F_\vt})u^*.
}\eeqn
\eitemp

Theorem \ref{strirepthm} can be specified further regarding uniqueness of $F_\vt$. First, recall that the inner product of $E$ generates a closed ideal $\cB_E:=\cls\AB{E,E}$ in $\cB$, the \hl{range ideal} of $E$. Further, we have $E^*\odot E=\cB_E$ by the canonical isomorphism $x^*\odot y\mapsto\AB{x,y}$. Finally, $\cB_E\odot E^*=E^*\odot E\odot E^*=E^*\odot\sK(E)=E^*$ so that $E^*$ may also be viewed as a correpondence from $\cB_E$ to $\sK(E)$. (The critical task, here, is nondegeneracy of the action of $\cB_E$.) It follows that $F_\vt$ may be viewed as a correspondence from $\cB_E$ to $\cC$. We collect some more results from \cite{MSS03p} in a theorem.

\bthm\label{unithm}
$F_\vt$ is the unique correspondence from $\cB$ to $\cC$ that is also a correspondence from $\cB_E$ to $\cC$. More precisely, if $\wt{F}$ is another correspondence from $\cB$ to $\cC$ and $\tilde{u}\colon E\odot\wt{F}\rightarrow F$ an isomorphism of correspondences from $\sB^a(E)$ to $\cC$, then
\beqn{
\AB{x_1,x_2}y
~\longmapsto~
x_1^*\odot\tilde{u}(x_2\odot y)
}\eeqn
establishes an isomorphism from the \nbd{\cB_E}\nbd{\cC}correspondence $\cls\cB_E\wt{F}=\cB_E\odot\wt{F}$ (tensor product over $\cB$, where $\cB_E$ is viewed as a correspondence from $\cB_E$ to $\cB$) to $F_\vt$.

In particular, if $E$ is \hl{full} (that is, if $\cB_E=\cB$), then $\wt{F}$ is isomorphic to $F_\vt$.
\ethm

\brem
Theorem \ref{strirepthm} in the case when $F$ is a Hilbert space and for representations of $\sF(E)$ rather than strict representations of $\sB^a(E)$ is Rieffel's \cite[Theorem 6.29]{Rie74}. What we added in \cite{MSS03p} is, essentially, the discussion for the extension to $\sB^a(E)$ and the generalization from Hilbert spaces to Hilbert modules for the representation space $F$. Modulo notation and technical discussion, the arguments in the proofs are the same.

The mechanism behind both proofs can be summarized by the observation that we may switch from modules over $\cB_E$ to modules over $\sK(E)$ and back by tensoring with $E$ or $E^*$ from the relevant side. In fact, the crucial identities $E^*\odot E=\cB_E$ and $E\odot E^*=\sK(E)$ mean just that $E$ is a Morita equivalence from $\sK(E)$ to $\cB_E$ and $E^*$ its inverse under tensor product. (A correspondence $E$ from $\cA$ to $\cB$, is a \hl{Morita equivalence} from $\cA$ to $\cB$, if $E^*\odot E=\cB$ and $E\odot E^*=\cA$. The latter means that the left action of $\cA$ on $E$ defines an isomorphism from $\cA$ onto $\sK(E)$. Tensoring a representation space $H$ of $\cB$ with a correspondence $E$ from $\cA$ to $\cB$ on the left gives a representation space $E\odot H$ of $\cA$ where the induced representation $a\odot\id_H$ acts. In fact, we obtain a functor from the category of representations of $\cB$ with intertwiners as morphisms to the corresponding category for $\cA$. This functor is an equivalence, if and only if $E$ is a Morita equivalence. As for \nbd{C^*}algebras not every equivalence is induced by a Morita equivalence, we treat, later, only the case of von Neumann algebras.)
\erem

Now we turn our interest to categories of Hilbert modules and functors between them. More precisely, we discuss the result that every sufficiently regular functor among them is implemented by tensoring from the right with a suitable fixed correspondence. The original version is Blecher's \it{Eilenberg-Watts theorem} \cite{Ble97}, proved there using operator space techniques. Here we use a different fomulation and discuss the approach from \cite{MSS03p} that does not involve operator spaces explicitly. (For experts: Blecher's categories have the same objects as ours but more morphisms. By definition, his functors restrict to ours and by an additional regularity condition, boundedness of the functor, that he must require, his functors are unique extensions of ours. So, essentially the two formulations coincide.)

Let $\cB$ be a \nbd{C^*}algebra. By $\eC_\cB^*$ we denote the category of Hilbert \nbd{\cB}modules with the adjointable mappings as morphisms. (Blecher would use the category $\eC_\cB$ with the same objects but bounded right linear mappings as morphisms.) A functor $\er\colon\eC_\cB^*\rightarrow\eC_\cC^*$ is a \hl{\nbd{*}functor} if $\er(a^*)=\er(a)^*$ for every morphism $a$. It is \hl{strict}, if the restriction $\er\upharpoonright\sB^a(E)$ is strict (in the sense explained above) for every object $E$. We wish to determine the structure of strict \nbd{*}functors. (Blecher would speak about \it{bounded} strict \nbd{*}functors, where \nbd{*}functor means in \cite{Ble97} that the restriction to our categories is a \nbd{*}functor in our sense.)

The key points of the approach in \cite{MSS03p} are as follows: By Theorem \ref{strirepthm} for every object $E$ there is a correspondence $F_E:=E^*\odot\er(E)$ and an isomorphism $u_E\colon E\odot F_E\rightarrow\er(F)$ such that $\er(a)=u_E(a\odot\id_{F_E})u_E^*$. Showing that for a direct sum $E_1\oplus E_2$, the submodule $u_{E_1\oplus E_2}(E_1\odot F_{E_1\oplus E_2})$ of $\er(E_1\oplus E_2)$ is canonically isomorphic to $u_{E_1}(E_1\odot F_{E_1})=\er(E_1)$ is a key point also in \cite{Ble97}. Thus, $E_1\odot F_{E_1\oplus E_2}\cong E_1\odot F_{E_1}$. In particular, by a (two-fold!) application of Theorem \ref{unithm} for an arbitrary object $E$, $F_\cB$ is the unique correspondence from $\cB$ to $\cC$ such that $\er(E)=u_E(E\odot F_E)\cong E\odot F_{E\oplus\cB}\cong E\odot F_\cB$. Therefore,  fixing the correspondence $F:=F_\cB$ from $\cB$ to $\cC$, we obtain that $\er(E)\cong E\odot F$. The following theorem, the \it{Eilenberg-Watts theorem}, fixes for each object a concrete isomorphism and states that the family of all these isomorphisms establishes a natural equivalence between the functors $\er$ and $\er_F:=\bullet\odot\id_F$.

By $_\cB\smash{\eC^*}\!\!\!_\cC$ we denote the category of correspondences from $\cB$ to $\cC$ with the bilinear adjointable mappings as morphisms.

\bitemp[Theorem \cite{MSS03p}.]\label{EWthm}
(\it{Eilenberg-Watts theorem}.) Let $\er\colon\eC^*_\cB\rightarrow\eC^*_\cC$ be a strict
\nbd{*}functor. Then $F=\cB^*\odot\er(\cB)$ is a correspondence in $_\cB\smash{\eC^*}\!\!\!_\cC$ such that the strict \nbd{*}functor $\er_F$, defined by setting $\er_F(E)=E\odot F$ and $\er(a)=a\odot\id_F$, is naturally equivalent to $\er$ via the \hl{natural transformation} given by the family of mappings $v_E\colon\er_F(E)\rightarrow\er(E)$ defined by setting $v_E(x\odot(b^*\odot z))=\er(xb^*)z$.

Moreover, $F$ is unique in $_\cB\smash{\eC^*}\!\!\!_\cC$. That is, if $\wt{F}\in{}_\cB\smash{\eC^*}\!\!\!_\cC$ is another correspondence such that $\er_{\wt{F}}$ is naturally equivalent to $\er$, then $\wt{F}\cong F$.
\eitemp

\brem
It is easy to show that the $v_E$ are isometries and that they fulfill the naturality condition $\er_F(a)=v_{E_2}^*\er(a)v_{E_1}$ ($a\in\sB^a(E_1,E_2)$). The discussion, as sketched before and detailed in \cite{MSS03p}, shows how to find the $v_E$ and that they are, indeed, surjective. See \cite{MSS03p} also for details that answer the question why Theorem \ref{EWthm} does not allow to conclude back easily to Theorem \ref{strirepthm}.
\erem

\brem
It is not difficult to show that, if $\vt_1$ and $\vt_2$ are strict homomorphisms that compose to $\vt_2\circ\vt_1$, then $F_{\vt_2\circ\vt_1}=F_{\vt_1}\odot F_{\vt_2}$. Also, the constructions iterate associatively. This leads, in particular, to the construction product systems of Hilbert modules (better: of correspondences) and was our original motivation to study homomorphisms $\vt$.
\erem

\bob
Theorems \ref{strirepthm}, \ref{unithm} and \ref{EWthm} have obvious generalizations to von Neumann modules, where $\vt$ and $\er$ will be normal, $F$ will be a von Neumann correspondence and tensor product are those in the category of von Neumann correspondences as discussed in the following section.
\eob

\section{Von Neumann modules and representations}\label{vNmsec}

Von Neumann algebras are algebras acting on a Hilbert space. It is easy to obtain them: Just take the closure of a \nbd{*}algebra of bounded operators or, equivalently, take its \it{double commutant}. The abstract (that is, without a defining representation by operators on a Hilbert space) counter part are \nbd{W^*}algebras. Once given such an algebra, there is no principal difference in difficulty between treating it as a von Neumann algebras or as a \nbd{W^*}algebras. (Although, we feel that the methods based on the theory of operators on Hilbert spaces with its topologies appears to be more direct than that what comes out if one tries to capture these topologies abstractly only in terms of the \nbd{W^*}algebra. But this is certainly only a matter of personal taste.) In this context, \it{the} basic result about \nbd{W^*}algebras is that given an abstract \nbd{*}algebra, if there is a possibility to turn it into a \nbd{W^*}algebra, then the way how to do it is unique. (\nbd{W^*}Algebras have unique pre-dual  Banach spaces.) The problems in the theory of \nbd{W^*}algebras occur when when we have \nbd{*}algebras that are not yet \nbd{W^*}algebras in this (unique) sense. It is, generally, a difficult task to find a good candidate for the (future) pre-dual Banach space (whose dual would, then, be the desired \nbd{W^*}algebra) \it{without} fixing a representation of the \nbd{*}algebra.

These problems pass directly over to \nbd{W^*}modules. Let $E$ be (pre-)Hilbert module over a \nbd{W^*}algebra $\cB$. Then $E$ is a \nbd{W^*}module, if it is \hl{self-dual} (that is, every bounded right linear mapping $E\rightarrow\cB$ arises as $x^*$ for a suitable $x\in E$) or, equivalently, if $E$ has a (unique) pre-dual Banach space. But, if $E$ is not yet self-dual (for instance, if $E$ is the Hilbert module tensor product of two \nbd{W^*}modules), then how to make it self-dual? Paschke \cite{Pas73} showed that every (pre-)Hilbert module over a \nbd{W^*}algebra admits a (unique) self-dual extension, but the explicit construction is not very handy. Rieffel \cite{Rie74a} showed how this extension can be obtained more easily, but only after fixing a faithful representation of the \nbd{W^*}algebra, that is, after having turned the \nbd{W^*}algebra into a von Neumann algebra. In fact, Rieffel's construction can be thought of as the beginning of the idea of a commutant for Hilbert modules.

So, if the (simple) construction of the (unique) self-dual extension depends on the choice of a representation, why not dealing from the beginning with von Neumann algebras rather than \nbd{W^*}algebras? Paired with the notion of von Neumann algebra as a concrete operator algebra, there is the notion of von Neumann module as a concrete module of operators. Much of the constructions needed to make this concept really applicable, can be found already in Rieffel \cite{Rie74a}, but the explicit definition (using strong closure in an operator space), that signifies a complete separation of the abstract properties from the concrete operator picture, is due to Skeide \cite{Ske00b}, and for the proof that this definition based on strong closure is equivalent to the one using self-duality we were not able to spot a reference going back further than \cite{Ske00b}.

\lf
So let $\cB\subset\sB(G)$ be a von Neumann algebra. (According to our convention, this means that $\cB$ is a strongly closed \nbd{*}algebra of bounded operators acting nondegenerately on the Hilbert space $G$.) We start by turning every (pre-)Hilbert \nbd{\cB}module $E$ into a concrete operator module. We define the Hilbert space $H=E\odot G$. Then every element $x\in E$ gives rise to an operator $L_x\colon g\mapsto x\odot g$ in $\sB(G,H)$ with adjoint $L_x^*$ defined by $y\odot g\mapsto\AB{x,y}g$. Clearly, $L_{xb}=L_xb$ and $\AB{x,y}=L_x^*L_y$. In other words, if we identify $x$ with $L_x$, then $E$ becomes a concrete operator \nbd{\cB}submodule of $\sB(G,H)$. We will always think in that way of $E$ as a subset of $\sB(G,H)$.

\bemp[Definition \cite{Ske00b}.]\label{vNmdef}
A (pre-)Hilbert module $E$ over a von Neumann algebra $\cB\subset\sB(G)$ is a \hl{von Neumann \nbd{\cB}module}, if $E$ is strongly closed in $\sB(G,H)$.
\eemp

\bcor\label{closcor}
If $E$ is a (pre-)Hilbert module over a von Neumann algebra $\cB\subset\sB(G)$, then the strong closure $\ol{E}^s$ in $\sB(G,H)$ is the unique smallest von Neumann \nbd{\cB}module containing $E$.
\ecor

\proof
This follows (like many other properties) simply because operator multiplication in $\sB(G\oplus H)$ is separately strongly continuous and $\cB$ is strongly closed in $\sB(G\oplus H)\supset\sB(G)$.\qed

\lf
Every $a\in\sB^a(E)$ gives rise to an operator in $\sB(H)$ that sends $x\odot g$ to $ax\odot g$. Instead of writing $a\odot\id_G$ we continue using the same letter $a$. In this way, we identify (faithfully, of course) $\sB^a(E)$ as a subalgebra of $\sB(H)$ acting nondegenerately on $H$. The following corollary follows as the preceding one.

\bcor
If $E$ is a von Neumann \nbd{\cB}module, then $\sB^a(E)\subset\sB(H)$ is a von Neumann algebra. (The converse need not be true, as the example $E=\sK(G,H)$ shows. In fact $\sK(G,H)$ is Hilbert \nbd{\sB(G)}module and $\sB^a(\sK(G,H))=\sB(H)$, but if $\sK(G,H)\ne\sB(G,H)$, that is, $H$ and $G$ are infinite-dimensional, then $\sK(G,H)$ is not a von Neumann module.)
\ecor

\bex
A von Neumann \nbd{\sB(G)}module has necessarily the form $E=\sB(G,H)$ (because it containes a subset norm-dense in the finite-rank operators $\sF(G,H)$) and $\sB^a(E)$ is $\sB(H)$. Therefore, who is interested in nontrivial operator algebras $\sB^a(E)\ncong\sB(H)$ may not look at \nbd{\sB(G)}modules.

Von Neumann modules $E$ over commutative von Neumann algebras (in standard representation, that is, $L^\infty$ acting by pointwise multiplication on $L^2$) as they occur in examples coming from classical probability, have direct integrals over the underlying measure space (type I von Neumann algebras) as operator algebras. In fact, Mingo and Giordano \cite{MiGio97} started considering this point of view as a possibility to free the theory of direct integrals from separability assumptions.
\eex

The basic result that makes the theory of von Neumann modules naturally equivalent to the theory of \nbd{W^*}modules is the following.

\bitemp[Theorem \cite{Ske00b,Ske03p}.]\label{vNsdthm}
A (pre-)Hilbert module over a von Neumann algebra is a von Neumann module, if and only if it is self-dual, that is, if and only if it is a \nbd{W^*}module.
\eitemp

The first proof in \cite{Ske00b} is based on existence of \it{quasi orthonormal bases}. The method of the second proof in \cite{Ske03p} is already closely related to the idea of commutant to which we gradually switch our attention.

We have imbedded $\sB^a(E)$ into $\sB(H)$ as acting on the first factor in $H=E\odot G$. This can be done for an arbitrary element in the algebra $\sB^r(E)$ of bounded right linear mappings on $E$.

\brem
The possibility to do so, is a nontrivial issue for the Banach algebra $\sB^r(E)$ without an \it{a priori} involution, and follows as, for instance, in Rieffel \cite{Rie74a} from general theorems about Banach modules. Only after this result, $\sB^r(E)$ turns out to be a not necessarily self-adjoint subalgebra of $\sB(H)$. The result is crucial for both proofs of the preceding theorem. In \cite{Ske00b} we provided a comparably elementary proof for existence and isometricity of the embedding $\sB^r(E)\subset\sB(H)$, that works only in the context of Hilbert modules and is based on \it{polar decomposition} and the \it{Kaplansky density theorem}.
\erem

It is natural to ask for the possibility to embed into $\sB(H)$ also operators that act on the second factor $G$ in $H=E\odot G$. This works, if the operator on $G$ is left \nbd{\cB}linear, that is, if it is an element of the commutant $\cB'$ of $\cB$. In other words, on $H$ we can define a unital (normal, of course) repesentation $\rho'$ of $\cB'$ by setting $\rho'(b')=\id_E\odot b'$ $(b'\in\cB')$. We call $\rho'$ the \hl{commutant lifting} associated with (pre-)Hilbert \nbd{\cB}module $E$. (This is not a one-hundred percent correct, looking at the meaning the term has in operator theory, but we think that \it{commutant lifting} expresses very well what $\rho'$ actually does; see also Example \ref{clStirem}.)

It is easy to compute the following commutant $\cM'$ and the double commutant $\cM''$ in $\sB(G\oplus H)=\sMatrix{\sB(G)&\sB(H,G)\\\sB(G,H)&\sB(H)}$ of the so-called \hl{linking algebra} $\cM=\sMatrix{\cB&E^*\\E&\sB^a(E)}$ (with obvious operations) of $E$:
\beq{\label{M'M''}
\cM'
~=~
\left\{\fMatrix{b'&0\\0&\rho'(b')}\colon b'\in\cB'\right\}
\text{~~~and~~~}
\cM''
~=~
\fMatrix{\cB&C_{\cB'}(\sB(H,G))\\C_{\cB'}(\sB(G,H))&\rho'(\cB')'},
}\eeq
where, generally, for an \nbd{\cA}bimodule $E$ we denote its \hl{\nbd{\cA}center} or just \hl{center} by
\beqn{
C_\cA(E)
~=~
\bCB{x\in E\colon ax=xa~(a\in\cA)}.
}\eeqn
So $C_{\cB'}(\sB(G,H))=\bCB{x\in\sB(G,H)\colon\rho'(b')x=xb'~(b'\in\cB')}$. As $\cM''$ is the strong closure of $\cM$, we conclude that $\ol{E}^s=C_{\cB'}(\sB(G,H))$ and and $\sB^a(\ol{E}^s)=\sB^a(E)''=\rho'(\cB')'$.

\bitemp[Corollary \cite{Ske03p}.]
$E$ is a von Neumann module, if and only if $E=C_{\cB'}(\sB(G,H))$. In this case, $\sB^a(E)=\rho'(\cB')'$.
\eitemp

Together with the result that $\sB^r(E)$ embeds into $\sB(H)$ it is not difficult to show that the bounded right linear mappings on $C_{\cB'}(\sB(G,H))$ embed into $C_{\cB'}(\sB(H,G))=C_{\cB'}(\sB(G,H))^*$, so that the intertwiner space $C_{\cB'}(\sB(G,H))$ is a self-dual Hilbert module. This was already observed by Rieffel \cite{Rie74a} and concludes the proof of self-duality of von Neumann modules ($E=C_{\cB'}(\sB(G,H))$!) as presented in \cite{Ske03p}.

\lf
Now a von Neumann \nbd{\cB}module gives rise to a representation of $\cB'$ and can be recovered as the intertwiner space for that representation. The question is natural, whether this correspondence can be reversed. That is, given a normal unital representation $\rho'$ of $\cB'$ on a Hilbert space $H$, can we define a von Neumann \nbd{\cB}module $E:=C_{\cB'}(\sB(G,H))$ so that its commutant lifting gives back $\rho'$. These are actually two questions. The first one, is $C_{\cB'}(\sB(G,H))$ a von Neumann \nbd{\cB}module, is readily verified to be affermative. (Excercise!) The second one, is $\rho'$ the commutant lifting associated with $E$, is tricky in two respects. Firstly, to construct the commutant lifting we have to construct $E\odot G$ and then $b'\mapsto\id_E\odot b'$. But, $E\odot G$ is a freshly constructed abstract space, while $H$ is given from the beginning. They cannot be \it{equal}, they can only be canonically isomorphic and $\rho'$ and the commutant lifting can, at most, be unitarilly equivalent. Secondly, suppose $E\odot G$ and $H$ are canonically isomorphic. Then, as $E$ generates all $E\odot G$ from $G$, the intertwiner space $C_{\cB'}(\sB(G,H))$ should do the same for $H$, that is, we should have $\cls C_{\cB'}(\sB(G,H))G=H$.

The first problem we resolve in a minute. (Under suitable specifications $H$ and $E\odot G$ are canonically isomorphic, and giving a suitable modified definition of \hl{concrete} von Neumann modules we sort this out.) The second problem has its affirmative solution in the following crucial lemma.

\bitemp[Lemma \cite{MuSo02}.]
If $\rho'$ is a normal unital representation of $\cB'$ on a Hilbert space $H$, then $\cls C_{\cB'}(\sB(G,H))G=H$.
\eitemp

Starting with $\rho'$ we define $\cM'$ as in \eqref{M'M''}. The idea of the proof is to see which closed subspace of $G\oplus H$ is generated from $G$ by the commutant $\cM''$ of $\cM'$ (as in \eqref{M'M''}). The projection $P'$ onto that subspace is in $\cM'$. So there exists a projection $p'\in\cB'$ such that $P'=\sMatrix{p'&0\\0&\rho'(p')}$. But $\cM''$ certainly generates all of $G$, so $p'=\U_{\cB'}$ and $P=\id_{G\oplus H}$.

For the solution of the first problem we collect the properties fulfilled by $E$ when identified as a subspace of $\sB(G,E\odot G)$, but formulate them in way where the Hilbert space $H$ is given from the beginning and $E$ is a concrete subset of $\sB(G,H)$.

\bdefi
Let $\cB\subset\sB(G)$ be a von Neumann algebra. A \hl{concrete von Neumann \nbd{\cB}mod\-ule} is a pair $(E,H)$ consisting of a Hilbert space $H$ and a subset $E$ of $\sB(G,H)$ such that:
\begin{enumerate}
\item\label{1}
$E$ is a (right) \nbd{\cB}submodule of $\sB(G,H)$, that is, $x\in E,b\in\cB$ $\Longrightarrow$ $xb\in E$.

\item\label{2}
$x,y\in E$ $\Longrightarrow$ $x^*y\in\cB$.

\item\label{3}
$E$ acts nondegeneratley on $G$, that is, $\cls EG=H$.

\item\label{4}
$E$ is strongly closed in $\sB(G,H)$.
\end{enumerate}
By $\cvN_\cB$ we denote the \hl{category of concrete von Neumann \nbd{\cB}modules} with the \hl{adjointable} mappings (that is a mapping $a\colon E_1\rightarrow E_2$ that admits a (unique) adjoint $a^*\colon E_2\rightarrow E_1$ such that $x^*(ay)=(a^*x)y$ for all $y\in E_1,x\in E_2$)  as morphisms.
\edefi

By \ref{1} and \ref{2} a concrete von Neumann \nbd{\cB}module $E$ is a pre-Hilbert \nbd{\cB}module with inner product defined as $\AB{x,y}=x^*y$.

By \ref{3} the Hilbert spaces $E\odot G$ and $H$ are isomorphic via the unitary defined by $x\odot g\mapsto xg$ (where \ref{3} contributes surjectivity). Therefore, every property present in the description using $H$ has its counterpart in the description using $E\odot G$. For instance, by \ref{4} the subset $E$ of $\sB(G,H)$ is strongly closed, thus, the same is true for the subset $\CB{L_x\colon x\in E}\subset\sB(G,E\odot G)$ so that $E$ is a von Neumann \nbd{\cB}module in the sense of Definition \ref{vNmdef}. It follows that the morphisms $E_1\rightarrow E_2$ are, indeed, the adjointable mappings $\sB^a(E_1,E_2)$ in the usual sense.

But, also each structure we defined so far in terms of $E\odot G$ has a counterpart when using $H$. The representation $b'\mapsto\id_E\odot b'$ of $\cB'$ on $E\odot G$ gives rise to a representation $\rho'\colon\cB'\rightarrow\sB(H)$ uniquely determined by $\rho'(b')xg=xb'g$. We recover $E$ as $E=C_{\cB'}(\sB(G,H))$. Moreover, the elements $a\in\sB^a(E_1,E_2)$ correspond one-to-one to elements in $C_{\cB'}(\sB(H_1,H_2))$, also denoted by $a$, where $a\in\sB^a(E_1,E_2)$ acts on $H_1$ as $a(x_1g)=(ax_1)g$ and where $a\in C_{\cB'}(\sB(H_1,H_2))$ acts on $x_1\in E_1\subset\sB(G,H_1)$ simply by composition $ax_1=a\circ x_1$ (but, except possibly in definitions, we never write the $\circ$ for compositions of operators on Hilbert spaces).

Conversely, if $(\rho',H)$ is a normal unital representation of $\cB'$ on a Hilbert space $H$, then $E:=C_\cB(\sB(G,H))$ defines a concrete von Neumann \nbd{\cB}module $(E,H)$ and the representation $b'\mapsto(xg\mapsto xb'g)$, constructed as before, gives us back $\rho'$. The correspondence between elements in $\sB^a(E_1,E_2)$ and in $C_{\cB'}(\sB(H_1,H_2))$ remains the same as discussed before. So, if we define $_{\cB'}\cvN$ as the \hl{category of normal unital representations of $\cB'$} with \nbd{\cB'}linear bounded (or, equivalently, adjointable) mappings as morphisms, then we obtain the following theorem.

\bthm\label{concthm}
Let $\cB\subset\sB(G)$ be a von Neumann algebra with commutant $\cB'\subset\sB(G)$. By
\baln{
\f
&
\colon
&
(E,H)
~&\longmapsto~
(b'\mapsto(xg\mapsto xb'g),H)
~,&~~~~~~~~
a
~&\longmapsto~
(a\colon x_1g\mapsto(ax_1)g)
\intertext{%
we define a bijective functor $\f\colon\cvN_\cB\rightarrow{}_{\cB'}\cvN$. The inverse functor is given by
}
\f^{-1}
&
\colon
&
(\rho',H)
~&\longmapsto~
(C_{\cB'}(\sB(G,H)),H)
~,&~~~~~~~~	
a
~&\longmapsto~
(a\colon x_1\mapsto a\circ x_1).
}\ealn
\ethm

We see that concrete von Neumann \nbd{\cB}modules and representations of $\cB'$ are isomorphic categories, not only naturally equivalent ones. Although the sequence $E\xrightarrow{~\f~}\rho'\xrightarrow{\f^{-1}}E$ certainly was known to Rieffel in \cite{Rie74a} and the back direction $\rho'\xrightarrow{\f^{-1}}E\xrightarrow{~\f~}\rho'$ must have been aware to many people working in Connes' setting of correspondences \cite{Con80p}, it seems that the one-to-one aspect of Theorem \ref{concthm}, featuring the usefulness of strongly closed operator modules as introduced (with emphasis on strong closure in the definition) in \cite{Ske00b}, has not been noticed so far.

\bob\label{modrepob}
Theorem \ref{concthm} tells us that speaking about von Neumann \nbd{\cB}modules and speaking about representations of $\cB'$ is the same thing. So, if $H$ is a Hilbert space, choosing a normal unital representation of $\cB'$ on $H$ turns it into a (concrete) von Neumann \nbd{\cB}module. In this picture, two von Neumann \nbd{\cB}modules determined by two representation $(\rho'_1,H_1)$ and $(\rho'_2,H_2)$ are established as isomorphic, if we, first, find a unitary $H_1\rightarrow H_2$ and, then, show that the unitary intertwines $\rho'_1$ and $\rho'_2$.

This harmless and obvious observation turns out to be a very powerful tool, when we have to identify von Neumann modules. The reason is, really, that it splits the construction of an isomorphism into two steps, the first of which is to well-define a mapping and, then, to show its properties. The definition of a mapping directly on the modules, usually, is somewhat in the converse order. First, one tries to give a prescription for how to calculate the mapping by phrasing the properties it should satisfy. In the first moment, one does not know whether the mapping is well-defined, but if it is, then it will have the desired properties. Only then, one shows that the mapping is well-defined.
\eob

We have learned that every statement or definition, for instance an \it{Eilenberg-Watts theorem}, for the category $\cvN_\cB$ can be translated into one for the other category $_{\cB'}\cvN$, and conversely, by conjugation with $\f$. We can already smell the relation between the two \it{Eilenberg-Watts theorems}, but both of them coinvolve a correspondence, namely, that which implements the respective functor. So, before we can really describe the relation, we have to discuss briefly von Neumann correspondences, in particular concrete ones, and we have to extend the pair of functors in Theorem \ref{concthm} to one functor, the commutant, sending (concrete) von Neumann \nbd{\cB}\nbd{\cC}correspondences to (concrete) von Neumann \nbd{\cC'}\nbd{\cB'}correspondences, and back. Anticipating the fact that a von Neumann \nbd{\cB}module is a von Neumann \nbd{\C}\nbd{\cB}correspondence, and that a representation of $\cB'$ is a von Neumann \nbd{\cB'}\nbd{\C'}correspondences ($\C'=\C\subset\sB(\C)=\C$), we will identify also $\f$ and $\f^{-1}$ as instances of the commutant.

\section{Von Neumann correspondences and their commutants}\label{commsec}

Let $\cA$ and $\cB$ denote \nbd{W^*}algebras. A correspondence $E$ from $\cA$ to $\cB$ is a \hl{\nbd{W^*}correspondence}, if $E$ is a \nbd{W^*}module over $\cB$ such that all the mappings $a\mapsto\AB{x,ax}$ $(x\in E)$ are normal. (Notice that we have nondegeneracy of the left action of $\cA$ in the \nbd{C^*}sense. But, as $\cA$ is unital, this does not matter.)

Let $\cA\subset\sB(K)$ and $\cB\subset\sB(G)$ denote von Neumann algebras. A correspondence $E$ from $\cA$ to $\cB$ is a \hl{von Neumann correspondence}, if $E$ is a von Neumann module over $\cB$ such that the canonical representation $\cA\rightarrow\sB^a(E)\subset\sB(H)$ is normal. We refer to $\rho\colon\cA\rightarrow\sB(H)$ as the \hl{Stinespring representation} of $\cA$ associated with $E$. It is routine to show that $E$ is a von Neumann correspondence, if and only if it is a \nbd{W^*}correspondence; see Skeide \cite[Lemma 3.3.2]{Ske01}.

\bex\label{clStirem}
Why are we refering to $\rho$ as the Stinespring representation? Because the construction of the original Stinespring representation \cite{Sti55} ``factors through Hilbert modules'' (or better through correspondences) and the way we defined $\rho$ captures exactly what happens. What do we mean by that?

Let $T\colon\cA\rightarrow\cB$ be a completely positive (CP-)map (for simlicity normal, between von Neuman algebras, but what we do works already for \nbd{C^*}algebras provided $\cB$ is represented faithfully on a Hilbert space $G$). Then Paschke's GNS-construction associates with $T$ a (unique) correspondence $E$ from $\cA$ to $\cB$ that is generated by a single vecctor $\xi\in E$ fulfilling $\AB{\xi,a\xi}=T(a)$ for all $a\in\cA$. ($E$ is simply the algebraic tensor product $\cA\otimes\cB$ with the only reasonable (semi-)inner product, length-zero elements quotiented out and completed.) Identifying $E\subset\sB(G,H)$, the representation $\rho$ is exactly Stinespring's representation on $H=E\odot G$ and $\xi=L_\xi\in\sB(G,H)$ the mapping such that $T(a)=\xi^*a\xi$.

And why do we refer to $\rho'$ as the commutant lifting? The strong closure of the GNS-correspondence in $\sB(G,H)$ is a von Neumann correspondence. And the representation $\rho'$ is exactly what Arveson is doing in the section with the title ``lifting commutants'' in \cite{Arv69}, when we interpret his CP-map $T\colon\cA\rightarrow\sB(G)$ as a mapping into $\cB:=T(\cA)''\subset\sB(G)$. (This is the minimal choice for $\cB$, consequently, with maximal commutant in $\sB(G)$. But, of course, we may choose for $\cB$ any von Neumann subalgebra of $\sB(G)$ that contains $T(\cA)$.)
\eex

The \hl{tensor product} $E\sodots F$ of von Neumann correspondences $E$ and $F$ is simply the strong closure of the usual tensor product $E\odot F$ in the sense of Corollary \ref{closcor}.  The corollary tells us that the strong closure is the unique self-adjoint extension and, therefore, coincides with the usual definition in the \nbd{W^*}framework. But, strong closure is much easier to obtain.

\bdefi
Let $\cA\subset\sB(K)$ and $\cB\subset\sB(G)$ be von Neumann algebras. A correspondence $E$ from $\cA$ to $\cB$ is a \hl{concrete von Neumann correspondence}, if $E$ is a concrete von Neumann module over $\cB$ and a von Neumann correspondence.

By $_\cA\cvN_\cB$ we denote the \hl{category of concrete von Neumann correspondences} from $\cA$ to $\cB$ with the bilinear adjointable mappings as morphisms. (Only left linearity must be checked, because right linearity follows from adjointability.)
\edefi

By Theorem \ref{concthm} a concrete von Neumann \nbd{\cB}module is given simply by a normal unital representation $(\rho',H)$, the commutant lifting. What we have to add to the representation $\rho'$ in order to have a concrete von Neumann correspondence from $\cA$ to $\cB$ is just a normal unital representation of $\cA$, the Stinespring representation. As $\rho$ maps into $\sB^a(E)=\rho'(\cB')'$, the two representations \hl{commute mutually}, that is, $\SB{\rho(\cA),\rho'(\cB')}=\zero$. But, this is the only condition a representation $\rho$ must satisfy in order to turn the concrete von Neumann module determined by $(\rho',H)$ into a concrete von Neumann correspondence.

\brem
A triple $(\rho',\rho,H)$ of a pair of normal unital mutually commuting representations $\rho'\colon\cB'\rightarrow\sB(H)$ and $\rho\colon\cA\rightarrow\sB(H)$ on a Hilbert space $H$ is very close to what Connes called \it{correspondence} in \cite{Con80p}. The missing link is as follows. If $\cB$ is in standard representation, then Tomita conjugation provides us with an isomorphism $\cB'\rightarrow\cB^{op}$. And a Connes correspondence from $\cA$ to $\cB$ is just a pair of commuting representations of $\cA$ and of $\cB^{op}$.

Our setting is slightly more general, but not much. (It is easy to show that two commutants of the same \nbd{W^*}algebra, obtained by choosing two faithful normal unital representations, are always Morita equivalent.) Appart from that, we think that our setting is considerably more elementary. We need not know what the standard representation is and we need not know (parts of) Tomita-Takesaki theory. We also mention that the construction of tensor products of such triples is a difficult task. Most discussions seem not to work without technical restrictions (typically to $\text{II}_1$ factors), while our definition of the tensor product is elementary, general and easily applicable.
\erem

Coming back to Theorem \ref{concthm} with the Stinespring representation of $\cA$ added, we have obtained a bijective functor from the category of concrete von Neumann \nbd{\cA}\nbd{\cB}correspondences $(E,H)$ to the category of triples $(\rho',\rho,H)$ with the mappings that intertwine both $\rho'$ and $\rho$ as morphisms.

So far, it was not really necessary to think of $\cA$ as a concrete von Neumann algebra acting on a Hilbert space $K$. However, as observed in Skeide \cite{Ske03c} (in the case $\cA=\cB$) and discussed also, independently, in Muhly and Solel \cite{MuSo04} in a \nbd{W^*}context (and generalized to different algebras in \cite{MuSo05}): In the triple picture $(\rho',\rho,H)$ the roles of $\rho'$ and of $\rho$ are in perfect symmetry. Nobody prevents us from switching to the triple $(\rho,\rho',H)$. Only when we go back to the correspondence picture the roles of $\cB'$ and $\cA$ are interchanged. Now $\cA$ is the commmutant of $\cA'$ and the representation $\rho$, interpreted as the commutant lifting of that commutant, determines a concrete von Neumann module $(E':=C_\cA(\sB(K,H)),K)$ over $\cA'$. Now it is the representation $\rho'$, when interpreted as Stinespring representation, that turns $E'$ into a correspondence from $\cB'$ to $\cA'$.

\bitemp[Definition and Theorem.]\label{commthm}
The diagram
\beqn{
\parbox{8cm}{
\xymatrix{
	&(\rho',\rho,H)	\ar@{<->}[dl]	\ar@3{-}[r]	&(\rho,\rho',H)	\ar@{<->}[dr]	&	\\
(E,H)	&									&									&(E',H)
}
}
}\eeqn
establishes a bijective functor $\f$, called the \hl{commutant functor}, from the category $_\cA\cvN_\cB$ of concrete von Neumann \nbd{\cA}\nbd{\cB}correspondences to the category $_{\cB'}\cvN_{\cA'}$ of concrete von Neumann \nbd{\cB'}\nbd{\cA'}correspondences.

We say $E'$ is the \hl{commutant} of $E$.

Varrying the parameters $\cA$ and $\cB$ of the functor to $\cB'$ and $\cA'$, we have $(E')'=E$.
\eitemp

\bob\label{corrrepob}
This is the counterpart of Observation \ref{modrepob}. Also correspondences, when thought of as triples $(\rho',\rho,H)$, can be identified by fixing first a unitary, and then showing that it intertwines now both representations.
\eob

If we identify $\cvN_\cB$ with $_\C\cvN_\cB$ in the only possible way and $_{\cB'}\cvN$ with $_{\cB'}\cvN_{\C}$ (identifying $H=\sB(\C,H)$ as $h\colon z\mapsto hz$), then the functors $\f$ and its inverse in Theorem \ref{concthm} become special cases of the commutant functor of Theorem \ref{commthm}. The commutant functor $\f$ is \it{cum grano salis} auto-inverse and we need no longer write $\f^{-1}$.

\section{Eilenberg-Watts theorems under commutant}\label{EWsec}

Let us describe the technical hypothesis of our functors in a unified way. We say that a functor $_\cA\cvN_\cB\rightarrow{}_\cC\cvN_\cD$ is \hl{normal} if for every object $E$ it restriction to $\sB^{a,bil}(E)$ is normal. A functor is a \hl{\nbd{*}functor}, if it respects adjoints.

The version for (concrete) von Neumann modules of Blecher's \it{Eilenberg-Watts theorem}, Theorem \ref{EWthm}, reads as follows.

\bthm\label{BEWthm}
Let $\cB\subset\sB(G)$ and $\cC\subset\sB(L)$ be von Neumann algebras and let $\er\colon\cvN_\cB\rightarrow\cvN_\cC$ be a normal \nbd{*}functor. Then there exists a unique up to isomorphism (concrete) von Neumann correpondence $F$ from $\cB$ to $\cC$ such that the functor
\beqn{
\er_F\colon
~~~~~~
E
~\longmapsto~
E\sbars{\odot}F,
~~~~~~
a
~\longmapsto~
a\odot\id_F
}\eeqn
is naturally equivalent to $\er$.
\ethm

Rieffel's \it{Eilenberg-Watts theorem}, in our language, takes the following form.

\bthm\label{REWthm}
Let $\cB\subset\sB(G)$ and $\cC\subset\sB(L)$ be von Neumann algebras and let $\el\colon{}_{\cB'}\cvN\rightarrow{}_{\cC'}\cvN$ be a normal \nbd{*}functor. Then there exists a unique up to isomorphism (concrete) von Neumann correpondence $F'$ from $\cC'$ to $\cB'$ such that the functor
\beqn{
\el_{F'}\colon
~~~~~~
H
~\longmapsto~
F'\odot H,
~~~~~~
a
~\longmapsto~
\id_{F'}\odot a
}\eeqn
is naturally equivalent to $\el$.
\ethm

Recall that $H$ is a Hilbert space so that $F'\odot H=F'\sodots H$. (Strong and norm topology on $\sB(\C,H)$ coincide. This is probably the reason, why a tensor product of von Neumann correspondences along the lines of Rieffel, that is, along the lines we described in the preceding section, has not been developed earlier.) Recall, too, that by Theorems \ref{concthm} and \ref{commthm} instead of $H$ we may write also $H=E'$ where $E=H'$. Then the two functors $\er_F$ and $\el_{F'}$ assume the perfectly symmetric form
\baln{
\er_F(E)
~&=~
E\sodots F
&
\el_{F'}(E')
~&=~
F'\sodots E'.
}\ealn
So far, $F$ and $F'$ are correspondences from different theorems. We just denoted the von Neumann \nbd{\cC'}\nbd{\cB'}correspondence granted by Theorem \ref{REWthm} by the symbol $F'$. As the $\er$ and $\el$ in the hypothesis are not related, $F$ from Theorem \ref{BEWthm} and $F'$ from Theorem \ref{REWthm} need not be related, either. But, suppose we have a functor $\er$, and we construct a functor $\el$ as $\el:=\f\circ\er\circ\f$, that is, $\el(E')=\er(E)'$. If we could show that the commutant functor takes tensor products to tensor products of the commutants in the opposite order (the only order that makes sense), that is, if we could show that
\beq{\label{tpcomm}
(E\sodots F)'
~\cong~
F'\sodots E',
}\eeq
then, indeed, $\el(E')=\er(E)'\cong\er_F(E)'=(E\sodots F)'\cong F'\sodots E'=\el_{F'}(E')$. As the (canonical) isomorphisms in this chain intertwine the actions of the relevant algebras, $\el(E')\cong\el_{F'}(E')$ provides us with a natural transform as claimed in Theorem \ref{REWthm}. In other words, Theorem \ref{REWthm} would follow from Theorem \ref{BEWthm}, but, as the whole discussion is symmetric, also Theorem \ref{BEWthm} would follow from Theorem \ref{REWthm}. We would, thus, have proved the following.

\bthm
Theorem \ref{BEWthm} and Theorem \ref{REWthm} are dual to each other under the commutant functor.
\ethm

Fortunately, Equation \ref{tpcomm} is true for arbitrary concrete von Neumann correspondences for which the tensor products make sense. This statement is \cite[Lemma 2.2]{Ske03c} or \cite[Lemma 3.7]{MuSo04} for a single von Neumann algebra (both based on a computation leading to \cite[Proposition 2.12]{MuSo02} but in their context there is no commutant of correspondences arround) and the general case in \cite[Lemma 3.9]{MuSo05}. The treatment in \cite{MuSo04,MuSo05} is rather on the level of \nbd{W^*}correspondences, while a treatment adapted exactly to the (concrete) von Neumann correspondences will appear in \cite{Ske03p1}. Here, in order to avoid messing up notation (by choosing many different repesentation spaces for many different von Neumann algebras) we treat only the case needed for \eqref{tpcomm} where one von Neumann algebra is $\C$. (But cf.\ also Remark \ref{funcomprem}.)

We would also like to mention that the definition of the tensor product of correspondences can be made such (choosing a different realization) that \eqref{tpcomm} becomes ``sharp'', that is, an equality and not just an isomorphism. But the discussion is tedious and we refer also here to \cite{Ske03p1}.

Let us come to the last missing piece.

\bthm\label{funflipthm}
Let $\cB\subset\sB(G)$ and $\cC\subset\sB(L)$ be von Neumann algebras. Then for every $(F,K)\in{}_\cB\cvN_\cC$ we have
\beqn{
\f\circ\er_F
~=~
\el_{F'}\circ\f
}\eeqn
up to natural equivalence.
\ethm

\proof
For every object $(E,H)\in\cvN_\cB$ we seek for identifications $\f\circ\er_F(E,H)=\el_{F'}\circ\f(E,H)$ by (canonical) isomorphisms (following Observations \ref{modrepob} and \ref{corrrepob}) that, then, provide automatically a natural transform. So, there is no harm if we decide to choose for the occuring tensor products a realization as concrete correspondences different from but isomorphic to that from the definition of the tensor product. Let $(\sigma',\sigma,K)$ denote the triple associated with $(F,K)$, that is, $\sigma'$ is the commutant lifting of $\cC'$ and $\sigma$ the Stinespring representation of $\cB$. For $E\sodots F$ we choose the realization as concrete von Neumann \nbd{\cC}module $(E\sodots F,E\odot K)$ in the following way. For every $x\in E$ define the mapping $\eta(x)\in\sB(K,E\odot K)$ (similar to $L_x$) by $k\mapsto x\odot k$. Observe that $E\odot K$ carries the representation $\tau'\colon c'\mapsto\id_E\odot\sigma'(c')$ of $\cC'$, and that $\eta(x)$ intertwines $\tau'$ and $\sigma'$. It follows that $E\sodots F:=\cls^s\eta(E)F\subset\sB(L,E\odot K)$ is a strongly closed subset of $C_{\cC'}(\sB(L,E\odot K)$ that has complement $\zero$ and, therefore, (like for every strongly closed submodule von Neumann module with zero-complement; see \cite[Corollary 3.2.12]{Ske01}) it follows that $E\sodots F=C_{\cC'}(\sB(L,E\odot K)$. On the other hand, the canonical isomorphism of $K$ and $F\odot L$ (in the sense of triples) shows that the von Neumann \nbd{\cC}module $E\sodots F$ is, indeed, isomorphic to the tensor product as defined in the previous section.

We find the equalities (without any canonical identification)
\beqn{
\f\circ\er_F(E,H)
~=~
\f(E\sodots F,E\odot K)
~=~
(\tau',E\odot K),
}\eeqn
and
\beqn{
\el_{F'}\circ\f(E,H)
~=~
\el_{F'}(\rho',H)
~=~
(\id_{\cC'}\odot\id_H,F'\odot H). 
}\eeqn
The spaces $E\odot K$ and $F'\odot H$ are different, even if we take into account that $K=\cls F'G$ and $H=\cls EG$. However, we have the canonical identifications $\cls F'G=F'\odot G$ and $\cls EG=E\odot G$. With these and the canonical identification $E\odot(F'\odot G)=F'\odot(E\odot G)$ defined by
\beq{\label{EE'tp}
x\odot y'\odot g
~\longmapsto~
y'\odot x\odot g
}\eeq
(Exercise: Check that \eqref{EE'tp} defines a unitary!) and taking also into account how $\cC'$ acts on these spaces, we find the desired identification up to canonical isomorphism. (Indeed, $c'$ on the left-hand side passes through $x$, because tensoring with $x$ may be replaced by the action of the intertwiner $\eta(x)$, so that $c'$ comes to act on $y'$ by left multiplication. Therefore, we obtain the same action as on the righ-hand side.)\qed

\brem
The operation described by \eqref{EE'tp} can be used to construct a tensor product of a von Neumann \nbd{\cB}module and a von Neumann \nbd{\cB'}module. The result is a von Neumann \nbd{(\cB\cap\cB')'}module, that is, a representation of the center $\cB\cap\cB'$ of $\cB$. This is closely related contructions on the Hilbert space level by Sauvageot \cite{Sau80,Sau83}. We describe details in \cite{Ske03p1}.
\erem

\brem\label{funcomprem}
Clearly, $\er_{F_2}\circ\er_{F_1}=\er_{F_1\sodots F_2}$ and $\el_{F'_1}\circ\el_{F'_2}=\el_{F'_1\sodots F'_2}$. Therefore, Theorem \ref{funflipthm} together with the uniquness of the correspondences inducing such functors may even be used to show in full generality that $(F_1\sodots F_2)'=F'_2\sodots F'_1$ up to isomorphism for arbitrary correspondences that match. But, in \cite{Ske03p1} we specify this better than just as up to isomorphism.
\erem

\brem
By uniqueness, the functor $\er_F$ ($\el_{F'}$) is an equivalence, if and only if $F$ ($F'$) is a Morita equivalence. This furnishes also new proofs for the corresponding \it{Morita theorems}.
\erem

\setlength{\baselineskip}{2.5ex}



\newcommand{\Swap}[2]{#2#1}\newcommand{\Sort}[1]{}
\providecommand{\bysame}{\leavevmode\hbox to3em{\hrulefill}\thinspace}
\providecommand{\MR}{\relax\ifhmode\unskip\space\fi MR }
\providecommand{\MRhref}[2]{%
  \href{http://www.ams.org/mathscinet-getitem?mr=#1}{#2}
}
\providecommand{\href}[2]{#2}


\end{document}